\documentclass[11pt,draft,amsmath,amscd,amsbsy,amssymb]{amsart}
 \relax

\theoremstyle{plain}
\newtheorem{thm}{Theorem}[section]
\newtheorem{cor}[thm]{Corollary}
\newtheorem{lemma}[thm]{Lemma}

\theoremstyle{definition}

\newcommand{\comp}{{\mathrm C}{\mathrm o}{\mathrm m}{\mathrm p} }
\newcommand{\tych}{{\mathrm T}{\mathrm y}{\mathrm c}{\mathrm h} }

\newcommand{\supp}{{\mathrm s}{\mathrm u}{\mathrm p}{\mathrm p} }
\newcommand{\card}{{\mathrm C}{\mathrm a}{\mathrm r}{\mathrm d} }
\newcommand{\pr}{{\mathrm p}{\mathrm r}}

\numberwithin{equation}{section}

\usepackage[all]{xy}

\begin{document}

\title[]
{On continuity of correspondences of probability measures in the
category of Tychonov spaces}
\author{Roman Kozhan}

\address{Department of Mechanics and Mathematics, Lviv National University,
Universytetska 1, 79000 Lviv, Ukraine} \email{}

\thanks{The author gratefully thanks Michael Zarichnyi for his
capable assistance in this research.}


\begin{abstract}We consider the correspondence assigning to
every Radon measure on two Tychonoff coordinate spaces the set of
probability measures with these marginals. It is proved that this
correspondence is continuous.
\end{abstract}

\maketitle

\section{Introduction}\label{s:intro}

In this paper we continue to investigate the relationship between
joint probability distribution on product spaces and their
marginal distributions. This problem is actual because of its
application in economic and game theories. In particular, we
consider the correspondence assigning to every probability measure
on two coordinate spaces the set of probability measures with
these marginals. Game theory regards a probability measure, for
instants, as a mixed strategy of players and a continuity of the
correspondence is a sufficient condition of existence of Nash
equilibrium of the static game in mixed strategies.

Similar problems were considered by Bergin (1999) for metrizable
coordinate spaces and Zarichnyi (2003) for the compact case. The
natural extension of these results is to take under consideration
spaces of probability measures in $\tych$. There are some
extensions of functor $P$ onto this category proposed by
Chigogidze (1984) and Fedorchuk (1991). In this paper we consider
the correspondence on the space of Radon measures and as a
corollary we derive the result for construction $P_{\beta}$.

The paper is organized as follows. Section \ref{s:defs} provides
necessary definitions and preliminary results which we need for
the proof of the main theorem. The continuity of correspondences
is established in Section \ref{s:mainres}.

\section{Definitions and preliminaries}\label{s:defs}

For the Tychonoff space $X$ let us consider the space
\[P_{\beta}(X)=\{\mu\in P(\beta X)\colon
\supp(\mu)\subset X\subset\beta X\},
\]where $\beta X$ is the Stone-Cech compactification  of $X$ and
$\supp(\mu)$ is the support of measure $\mu$. This construction
extends functor $P\colon\comp\to\comp$ to the functor
$P_{\beta}\colon\tych\to\tych$ and was proposed by Chigogidze
(1984) (see also Teleiko and Zarichnyi (1999)).

In order to consider two other extensions let us recall some
definitions. The set of all Borel measures on Tychonoff space $X$
is denoted by $M(X)$. The weak-* topology on $M(X)$ is the
topology generated by the base which consists of the sets of the
form
\[O(\mu, V_1,...,V_k,
\varepsilon)=\{\nu\in\hat{P}(X)\colon\nu(V_i)>\mu(V_i)-\varepsilon,
i=1,...,k\},
\]where $V_i\subset X$ are open sets and $\varepsilon>0$.

The Borel measure $\mu\in M(X)$ is called Radon measure if
\[\mu(A)=\sup\{\mu(K)\colon A\supset K -\text{ compact subset of
X}\}
\]for every Borel subset $A\subset X$.

The Borel measure $\mu\in M(X)$ is called $\tau$-smooth if for
every decreasing sequence $\{Z_{\alpha}\}$ of closed subsets of
$X$ such that $\underset{\alpha}{\cap}Z_{\alpha}=\emptyset$ the
sequence $\{\mu(Z_{\alpha})\}$ approaches to zero.

The subset of $M(X)$ formed by all Radon measures is denoted by
$\hat{M}(X)$ and $M_{\tau}(X)$ defines the set of all
$\tau$-smooth Borel measures. We denote by $\hat{P}(X)$ and
$P_{\tau}(X)$ subsets of Radon and $\tau$-smooth probability
measures of $X$ respectively. For topological properties of the
spaces $\hat{P}(X)$ and $P_{\tau}(X)$ see Banakh (1995).

Let $X$ and $Y$ be Tychonoff spaces.

By Riez representation theorem there is one-to-one map between
measures on the topological space $X$ and linear functionals on
the space of all bounded continuous functions $C_b(X)$ equipped
with the sup-norm. Varadarajan (1961) defines analogically
$\tau$-smooth and dense functionals.

Recall, that bounded linear functional $\Lambda$ is called
$\tau$-smooth if $f_{\alpha}\downarrow 0$ implies that
$\Lambda(f_{\alpha})\to 0$ for every net $\{f_{\alpha}\}\subset
C_b(X)$.

Similarly, bounded linear functional $\Lambda$ is called dense if
$\Lambda(f_{\alpha})\to 0$ for any net $\{f_{\alpha}\}\subset
C_b(X)$ such that $\|f_{\alpha}\|\leq 1$ for every $\alpha$ and
$f_{\alpha}\to 0$ uniformly on compact subsets of $X$.

Varadarajan (1961) proved that the measure $\mu$ is Radon
($\tau$-smooth) if and only if there exists dense ($\tau$-smooth)
bounded linear functional $\Lambda$ such that for every $f\in
C_b(X)$ we have $\Lambda(f)=\underset{X}{\int}fd\mu$ and
$\|\Lambda\|=\|\mu\|$.

Further, we will assume these notions equivalent and use both
techniques without comments.

Denote by $u_X\colon C(\beta X)\to C(P(\beta X))$ the map defined
by the formula $u_X(\varphi)(\mu)=\mu(\varphi)$. The map is a
linear operator with $\|u_X\|=1$.

For every $\varphi\in C_b(X)$ we can restrict the function
$u_X(\varphi)$ on $\hat{P}(X)$.

Define the map $\psi_X \colon P^2(\beta X)\to P(\beta X)$ by the
formula $\psi_X(M)(\varphi)=M(u_X(\varphi))$.

\begin{lemma}\label{l:1} $\psi_X(\hat{P}^2(X))\subseteq
\hat{P}(X)$.
\end{lemma}

\begin{proof} Consider the sequence of functions
$\{\varphi_\alpha\}_{\alpha\in \Gamma}$ such that
$\|\varphi_\alpha\|\leq 1$ for every $\alpha\in \Gamma$ and
$\varphi_\alpha \rightarrow 0$ uniformly on compact subsets of
$X$. It is easy to show that $\|u_X (\varphi_\alpha)\|\leq 1$ for
every $\alpha\in \Gamma$ and $u_X(\varphi_\alpha)\to 0$ uniformly
on compact subsets of $\hat{P}(X)$. Therefore for every radon
measure $M\in \hat{P}^2(X)$ we have
\[\psi_X(M)(\varphi_\alpha)=M(u_X(\varphi_\alpha))\to 0
\]and this implies that $\psi_X(M)\in\hat{P}(X)$.
\end{proof}

\begin{lemma}\label{l:2} i) For every $\mu\in
\hat{M}(X)$ and $\nu\in \hat{M}(Y)$ such that $\|\mu\|=\|\nu\|$
there exists $\lambda\in \hat{M}(X\times Y)$ such that
$\hat{M}\pr_1(\lambda)=\mu$ and $\hat{M}\pr_2(\lambda)=\nu$.

ii) The statement is true if we substitute $\hat{M}$ with
$M_{\beta}$.
\end{lemma}

\begin {proof} i). Let us first assume that
$\mu\in\hat{P}(X)$ and $\nu\in\hat{P}(Y)$. Denote by $i_x\colon
Y\to X\times Y$ an embedding defined for every $x\in X$ by the
formula $i_x(y)=(x,y)$. Define $f_\nu\colon X\to \hat{P}(X\times
Y)$ in the following manner: $f_\nu(x)=\hat{P}i_x(\nu)$ for every
$x\in X$. Let us check whether the measure $\lambda=\psi_{X\times
Y}(\hat{P}f_\nu(\mu))\in \hat{P}(X\times Y)$ satisfies the
conditions of the lemma. For every $f\in C_b(X)$
\begin{eqnarray*}
\hat{P}\pr_1(\lambda)(f)&=&\psi_{X\times
Y}(\hat{P}f_\nu(\mu))(f\circ \pr_1)=\hat{P}f_\nu(\mu)(u_{x\times
Y}(f\circ \pr_1))\\&=&\mu(u_{X\times Y}(f\circ\pr_1)\circ f_\nu).
\end{eqnarray*}Since the function
\begin{eqnarray*}
(u_{X\times Y}(f\circ\pr_1)\circ f_\nu)(x)&=&u_{X\times Y}(f\circ
\pr_1)(f_\nu)=u_{X\times Y}(f\circ\pr_1)( \hat{P}i_x(\nu))\\
&=&\hat{P}i_x(\nu)(f\circ\pr_1)=\nu(f\circ\pr_1\circ i_x)
=\nu(f(x))=f(x)
\end{eqnarray*}we have that $u_{X\times Y}(f\circ\pr_1)\circ
f_\nu\equiv f$ and $\hat{P}\pr_1(\lambda)=\mu$.

Analogically for every $g\in C_b(Y)$
\begin{eqnarray*}
\hat{P}\pr_2(\lambda)(g)&=&\psi_{X\times
Y}(\hat{P}f_\nu(\mu))(f\circ\pr_1)=\mu(u_{X\times Y}(g\circ
\pr_2)\circ f_\nu).
\end{eqnarray*}Since the function
\begin{eqnarray*}
(u_{X\times Y}(g\circ\pr_2)\circ f_\nu)(x)&=&\nu(g\circ\pr_2\circ
i_x)=\nu(g)
\end{eqnarray*}this implies that $\hat{P}\pr_2(\lambda)(g)=
\mu(\nu(g))=\nu(g)$.

According to Lemma \ref{l:1} $\lambda\in \hat{P}(X\times Y)$. We
denote this measure by $\lambda=\mu\otimes\nu$ and call the tensor
product of the measures $\mu$ and $\nu$.

Now, if $c=\|\mu\|=\|\nu\|\neq 1$ then it is easy to check that
the measure $c(\frac{\mu}{c}\otimes\frac{\nu}{c})$ satisfies the
conditions of the lemma.

ii). Let us denote by $K_\mu=\supp(\mu)$ and $K_\nu=\supp(\nu)$.
By the definition of the functor $P_\beta$ the sets $K_\mu$ and
$K_\nu$ are compacta. This implies that measures $\frac{\mu}{c}\in
P(K_\mu)$ and $\frac{\nu}{c}\in P(K_\nu)$. The statement ii) of
the lemma follows from the bicommutativity of the functor $P$.
\end{proof}

\begin{lemma}\label{l:3} Let $\lambda\in \hat{P}(X\times Y)$.
Then for every open set $V\subset X\times Y$ and
$\varepsilon>0$ there exist elements of the base
$\tilde{V}_1,...,\tilde{V}_k$ such that $\tilde{V}_i\subset V$ for
every $i=1,...,k$ and
\[\lambda(V)-\lambda(\underset{i=1}{\overset{k}{\cup}}
\tilde{V}_i)<\varepsilon.
\]
\end{lemma}

\begin{proof} The sets $\tilde{V}_1,...,\tilde{V}_k$ are
elements of the base means that there are open sets
$\tilde{V}^{\prime}_i \subset X$ and
$\tilde{V}^{\prime\prime}_i\subset Y$ such that
$\tilde{V}_i=\tilde{V}^{\prime}_i\times \tilde{V}^{\prime\prime}_i
$ for every $i=1,...,k$. Since $\lambda$ is a Radon measure there
exists a compactum $K\subset V$ such that $\lambda(V\setminus
K)<\varepsilon$. Consider a covering by elements of the base
$\{\tilde{V}_i\}_{i\in\Gamma}$ of the compactum $K$ such that
$\tilde{V}_i\subset V$ for every $i\in\Gamma$. We can find a
finite number of elements of the covering $\tilde{V}_1,...,
\tilde{V}_k$ such that
$K\subset\underset{i=1}{\overset{k}{\cup}}\tilde{V}_i\subset V$.
This condition proves the lemma.
\end{proof}

\begin{lemma}\label{l:4} Let $\lambda\in \hat{P}(X\times
Y)$, $\mu=\hat{P}\pr_1(\lambda)$, $\nu=\hat{P}\pr_2(\lambda)$ and
$\varepsilon>0$. Then for open sets $V,V^{\prime}\subset X$ such
that $\mu(V\setminus V^{\prime})<\varepsilon$ it is satisfied that
$\lambda((V\setminus V^{\prime})\times W)<\varepsilon$ for every
open set $W\subset Y$. Analogically, for open sets
$W,W^{\prime}\subset Y$ such that $\nu(W\setminus
W^{\prime})<\varepsilon$ it is satisfied that
$\lambda(V\times(W\setminus W^{\prime}))<\varepsilon$ for every
open set $V\subset Y$.
\end{lemma}

\begin{proof} Since $\mu=\hat{P}\pr_1(\lambda)$
this implies that $\mu(V)=\lambda(\pr_1^{-1}(V))=\lambda(V\times
Y)$. Then for every open set $W\subset Y$ we have that
\[\lambda((V\setminus V^{\prime})\times
W)<\lambda((V\setminus V^{\prime})\times Y)=\mu(V\setminus
V^{\prime})<\varepsilon.
\] Another statement of the lemma can be proved in the same
manner.
\end{proof}

\begin{lemma}\label{l:5} $\lambda\in \hat{P}(X\times Y)$,
$\mu=\hat{P}\pr_1(\lambda)$, $\nu=\hat{P}\pr_2(\lambda)$,
$\varepsilon_1>0$, $\varepsilon_2>0$. Then for open sets $V,
V^{\prime}\subset X$ such that $\mu(V\setminus
V^{\prime})<\varepsilon_1$ and for open sets $W,W^{\prime}\subset
Y$ such that $\nu(W\setminus W^{\prime})<\varepsilon_2$ it is
satisfied that $\lambda((V\times W)\setminus(V^{\prime}\times
W^{\prime}))<\varepsilon_1+\varepsilon_2$.
\end{lemma}

\begin{proof} It is easy to check that
\begin{eqnarray*}
\lambda((V\times W)\setminus(V^{\prime}\times W^{\prime}))
&<&\lambda(((V\setminus V^{\prime})\times
W)\cup(V\times(W\setminus W^{\prime})))\\
&\leq&\lambda((V\setminus V^{\prime})\times
W)+\lambda(V\times(W\setminus W^{\prime}))<
\varepsilon_1+\varepsilon_2.
\end{eqnarray*}
\end{proof}

\begin{lemma}\label{l:6} Let $\lambda\in\hat{P}(X\times
Y)$ and let $V_1,...,V_n$ be open sets in $X\times Y$ and
$\varepsilon_0>0$. Then for every $\varepsilon>0$ there exist
pairwise disjoint elements of base $W_1,...,W_m\subset X\times Y$
and a number $\delta>0$ such that

(i) $\lambda(V_i)-\underset{\{j\colon W_j\subset
V_i\}}{\sum}\lambda(W_j)<\varepsilon$ for every $i=1,...,n,$

(ii) $O(\lambda,W_1,...,W_m,\delta)\subset
O(\lambda,V_1,...,V_n,\varepsilon_0)$,

(iii) $(\pr_l(W_{j^{\prime}})\cap \pr_l(W_{j^{\prime\prime}})
=\emptyset)\vee (\pr_l(W_{j^{\prime}})=
\pr_l(W_{j^{\prime\prime}}))$ for $l=1,2$ and
$j^{\prime},j^{\prime\prime}=1,...,m$.
\end{lemma}

\begin{proof} According to Lemma 2.2 of Bogachev (1999) we can
assume that $V_1,...,V_n$ are pairwise disjoint. By Lemma
\ref{l:2} we can find elements of base $\tilde{V_{ij}}\subset
V_i$, $j=1,...,k_i$ such that
\[\lambda(V_i)-\lambda(\underset{j=1}{\overset{k_i}{\cup}}\tilde{V}_{ij})
<\frac{\varepsilon}{2}.
\]Since every set $\tilde{V}_{ij}$ is an
element of the base it can be represented as
$\tilde{V}_{ij}=\tilde{V}_{ij}^{\prime}\times\tilde{V}_{ij}^{\prime\prime}$
for some open sets $\tilde{V}_{ij}^{\prime}\subset X$ and $
\tilde{V}_{ij}^{\prime\prime}\subset Y$. Following Lemma 2.2 of
Bogachev (1999) there exist pairwise disjoint open sets
$W_1^{\prime},...,W_{m^{\prime}}^{\prime}\subset X$ and
$W_1^{\prime\prime},...,W_{m^{\prime\prime}}^{\prime\prime}\subset
Y$ such that
\[\mu(\tilde{V}_{ij}^{\prime})-\underset{\{s\colon
W_s^{\prime}\subset\tilde{V}_{ij}^{\prime}\}}{\sum}\mu(W_s^{\prime})<\frac{\varepsilon}{4k}
\]and
\[\nu(\tilde{V}_{ij}^{\prime\prime})-\underset{\{t\colon
W_t^{\prime\prime}\subset\tilde{V}_{ij}^{\prime\prime}\}}{\sum}
\nu(W_t^{\prime\prime})<\frac{\varepsilon}{4k}
\]for every $i=1,...,n$ with $k=\max\{k_1,...,k_n\}$. Consider all pairwise products of the type
$W_q^{\prime}\times W_s^{\prime\prime}$ for $q=1,...,m^{\prime}$,
$s=1,...,m^{\prime\prime}$ which is in at least one of the set
$\tilde{V}_{ij}$ for $j=1,...,k_i$ and $i=1,...,n$. Let us denote
these products as $W_1,...,W_m$. Clearly, they are pairwise
disjoint and $I_j=\{i\colon W_j\subset V_i\}$ is a singleton. In
addition, Lemma \ref{l:5} implies that for every $j=1,...,k_i$
\begin{eqnarray*}
\lambda(\tilde{V}_{ij})-\underset{\{t\colon
W_t\subset\tilde{V}_{ij}\}}{\sum}\lambda(W_t)&=&
\lambda(\tilde{V}_{ij}^{\prime}\times\tilde{V}_{ij}^{\prime\prime}\setminus
((\underset{\{q\colon W_q^{\prime}\subset\tilde{V}_{ij}^{\prime}
\}}{\cup}W_q^{\prime})\times(\underset{\{s\colon
W_s^{\prime\prime}\subset\tilde{V}_{ij}^{\prime\prime}
\}}{\cup}W_s^{\prime\prime})))\\&<&\frac{\varepsilon}{4k}+
\frac{\varepsilon}{4k}<\frac{\varepsilon}{2k}.
\end{eqnarray*}Let us show that for every $i=1,...,n$ we have
$\lambda(\underset{j=1}{\overset{k_i}{\cup}}\tilde{V}_{ij})-
\underset{q=1}{\overset{m}{\sum}}\lambda(W_q)<\varepsilon$. For
$j=1,...,k_i$ the inequality
\[\lambda(\tilde{V}_{ij}\setminus(\underset{q=1}{\overset{m}{\cup}}W_q))
\leq\lambda(\tilde{V}_{ij}\setminus(\underset{\{q\colon
W_q\subset\tilde{V}_{ij}\}}{\cup}W_q))<\frac{\varepsilon}{2k}.
\]is satisfied. Using this fact we can see that
\[\lambda((\underset{j=1}{\overset{k_i}{\cup}}\tilde{V}_{ij})\setminus
(\underset{q=1}{\overset{m}{\cup}}W_q))\leq
\underset{j=1}{\overset{k_i}{\sum}}\lambda(\tilde{V}_{ij}\setminus
(\underset{q=1}{\overset{m}{\cup}}W_q))<k_i\frac{\varepsilon}{2k}<\frac{\varepsilon}{2}.
\]The condition (i) is easily implied from the last equation.

Indeed, for every $i=1,...,n$
\[\lambda(V_i)-\underset{q=1}{\overset{m}{\sum}}\lambda(W_q)=
\lambda(V_i)-\lambda(\underset{j=1}{\overset{k_i}{\cup}}\tilde{V}_{ij})+
\lambda(\underset{j=1}{\overset{k_i}{\cup}}\tilde{V}_{ij})-
\underset{q=1}{\overset{m}{\sum}}\lambda(W_q)<\frac{\varepsilon}{2}+\frac{\varepsilon}{2}=
\varepsilon.
\]

In order to prove (ii) it is necessary to find $\delta$ such that
for $i=1,...,n$
\[\lambda^{\prime}(V_i)-\lambda(V_i)>-\varepsilon_0
\]for every $\lambda^{\prime}\in O(\lambda,W_1,...,W_m,\delta)$. The
condition (i) implies that for $\delta<\frac{\varepsilon_0}{2m}$
and $\varepsilon<\frac{\varepsilon_0}{2}$ we have
\begin{eqnarray*}\lambda^{\prime}(V_i)-\lambda(V_i)&>&\underset{\{q\colon
W_q\in V_i\}}{\sum}\lambda^{\prime}(W_q)-
\underset{\{q\colon W_q\in V_i\}}{\sum}\lambda(W_q)-\varepsilon\\
&>&\underset{\{q\colon W_q\in V_i\}}{\sum}
(\lambda^{\prime}(W_q)-\lambda(W_q))-\varepsilon>
-m\delta-\varepsilon>-\varepsilon_0.
\end{eqnarray*}

The statement (iii) is obvious by the algorithm.

The Lemma is proved.
\end{proof}

Further we need the notion of the restriction of the measure. We
endow the set
\[\{\psi\in C_b(X,[0,1]), \psi|_A\equiv 1\}
\]with
the reverse pointwise order relation. For every $\mu\in M(X)$ and
Borel subset $A\subset X$ define
$\mu|_A(\varphi)=\lim\{\mu(\varphi\cdot\psi\colon\psi\in
C_b(X,[0,1])), \psi|_A\equiv 1\}$ for every $\varphi\in C_b(X)$
(see Teleiko and Zarichnyi (1999)). It is easy to show that for
every Borel set $B\subset X$ we have that $\mu|_A(B)=\mu(A\cap
B)$.

Let us define the map $\hat{\chi}\colon \hat{P}(X\times Y)\to
\hat{P}(X)\times\hat{P}(Y)$ by the formula
\[\hat{\chi}(\lambda)=(\hat{P}\pr_1(\lambda),
\hat{P}\pr_2(\lambda))
\] for every
$\lambda\in\hat{P}(X\times Y)$.

\section{Main results}\label{s:mainres}

\begin{thm}\label{t:maintych} The map $\hat{\chi}$ is open.
\end{thm}

\begin{proof} Let $\lambda^0\in\hat{P}(X\times Y)$
and let $O(\lambda^0,V_1,...,V_n,\varepsilon)$ be its neighborhood
in the weak-* topology for open sets $V_i\subset X\times Y$ and
$\varepsilon>0$. According to Lemma \ref{l:6} we assume without
loss of generality that the sets $V_i$ are pairwise disjoint and
there are exist open sets $W_i^{\prime}\subset X$ and
$W_i^{\prime\prime}\subset Y$ such that $V_i=W_i^{\prime}\times
W_i^{\prime\prime}$ and $( W_{i_1}^{\prime}\cap
W_{i_2}^{\prime}=\emptyset)\vee
(W_{i_1}^{\prime}=W_{i_2}^{\prime})$ and $( W_{i_1}^{\prime\prime}
\cap W_{i_2}^{\prime\prime}=\emptyset)
\vee(W_{i_1}^{\prime\prime}=W_{i_2}^{\prime\prime})$ for
$i,i_1,i_2=1,...,n$. Let $\mu^0=\hat{P}\pr_1(\lambda^0)$ and
$\nu^0=\hat{P}\pr_2(\lambda^0)$. Denote by
$m^{\prime}=\card\{W_i^{\prime},i=1,...,n\}$ and
$m^{\prime\prime}=\card\{W_i^{\prime\prime},i=1,...,n\}$ and let
$\delta<\varepsilon$. We adopt the notation
$W_{qs}=W^{\prime}_q\times W^{\prime\prime}_s$ for every
$q\in\{1,...,m^{\prime}\}$ and $s\in\{1,...,m^{\prime\prime}\}$.

In order to prove the theorem it is sufficient to show that every
pair of measures
\[(\mu,\nu)\in O(\mu^0,W_1^{\prime},...,
W_{m^{\prime}}^{\prime},\delta)\times
O(\nu^0,W_1^{\prime\prime},...,W_{m^{\prime\prime}}^{\prime\prime},\delta)
\]has a preimage in $O(\lambda^0,V_1,...,V_n,\varepsilon)$.

For every
$(q,s)\in\{1,...,m^{\prime}\}\times\{1,...,m^{\prime\prime}\}$ we
define
\[\alpha^{\prime}_{qs}=\frac{\lambda^0(W_{qs})\mu(W_q^{\prime})}{
\mu^0(W^{\prime}_q)}
\]and analogically
\[\alpha^{\prime\prime}_{qs}=\frac{\lambda^0(W_{qs})\nu(W_s^{\prime\prime})}{
\nu^0(W^{\prime\prime}_s)}.
\]In addition we set
$\alpha_{qs}=\min(\alpha_{qs}^{\prime},\alpha_{qs}^{\prime\prime})$,
$\beta_{qs}^{\prime}=\alpha_{qs}^{\prime}-\alpha_{qs}$ and
$\beta_{qs}^{\prime\prime}=\alpha_{qs}^{\prime\prime}-\alpha_{qs}$.

Let $\mu_q=\mu|_{W^{\prime}_q}$ and
$\nu_s=\nu|_{W^{\prime\prime}_s}$ be the restrictions of the
measures $\mu$ and $\nu$ on the sets $W^{\prime}_q$ and
$W^{\prime\prime}_s$ respectively for every
$q\in\{1,...,m^{\prime}\}$ and $s\in\{1,...,m^{\prime\prime}\}$.
Denote
\[\lambda_{qs}=\alpha_{qs}(\frac{\mu_q}
{\mu(W^{\prime}_q)}\otimes \frac{\nu_s}{\mu(W^{\prime\prime}_s)})
\]if $\lambda^0(W^{\prime}_q\times W^{\prime\prime}_s)\neq
0$ and $\lambda_{qs}=0$ otherwise. It is clear that
\[\lambda_{qs}(X\times Y)=\lambda_{qs}(W^{\prime}_q \times
W^{\prime\prime}_s)=\alpha_{qs}.
\]

We define the measures
\[\tilde{\lambda}=\underset{q=1}{\overset{m^{\prime}}{\sum}}
\underset{s=1}{\overset{m^{\prime\prime}}{\sum}}\lambda_{qs}\in\hat{M}(X\times
Y)\text{, }\tilde{\mu}=\mu-\hat{M}\pr_1(\tilde{\lambda})
\in\hat{M}(X)\text{ and }\tilde{\nu}=\nu-\hat{M}\pr_2
(\tilde{\lambda})\in\hat{M}(Y).
\]Since for every $q\in\{1,...,m^{\prime}\}$ we have
\begin{eqnarray*}\tilde{\mu}(W^{\prime}_q)&=&\mu(W^{\prime}_q)-
\underset{s=1}{\overset{m^{\prime\prime}}{\sum}}\alpha_{qs}\geq
\mu(W^{\prime}_q)-\underset{s=1}{\overset{m^{\prime\prime}}{\sum}}
\alpha^{\prime}_{qs}\\&=&
\frac{\mu(W^{\prime}_q)\lambda^0(W^{\prime}_q\times
(Y\setminus(\underset{s=1}{\overset{m^{\prime\prime}}{\cup}}
W^{\prime\prime}_s)))}{\mu^0(W^{\prime}_q)}\geq 0.
\end{eqnarray*}this implies that the measures $\tilde{\mu}$ and $\tilde{\nu}$
are non-negative. Moreover, it is clear that
$\|\tilde{\mu}\|=\|\tilde{\nu}\|\geq 0$. Therefore by Lemma
\ref{l:2} there exists non-negative measure
$\tilde{\tilde{\lambda}}$ such that
$\hat{M}\pr_1(\tilde{\tilde{\lambda}})=\tilde{\mu}$ and
$\hat{M}\pr_2(\tilde{\tilde{\lambda}})=\tilde{\nu}$. Let us set
$\lambda=\tilde{\lambda}+\tilde{\tilde{\lambda}}$. Obviously,
$\lambda\in\hat{P}(X\times Y)$. Besides, for every
$q\in\{1,...,m^{\prime}\}$ and $s\in\{1,...,m^{\prime\prime}\}$ we
have
\begin{eqnarray*}\lambda(W_{qs})-\lambda^0(W_{qs})&=&
\tilde{\lambda}(W_{qs})+\tilde{\tilde{\lambda}}(W_{qs})
-\lambda^0(W_{qs})\geq \tilde{\lambda}(W_{qs})-\lambda^0(W_{qs})
\\&=&\alpha_{qs}-\lambda^0(W_{qs})=\lambda^0(W_{qs})(
\min\{\frac{\mu(W^{\prime}_q)}{\mu^0(W^{\prime}_q)},
\frac{\nu(W^{\prime\prime}_s)}{\nu^0(W^{\prime\prime}_s)} \}-1).
\end{eqnarray*}The
condition $(\mu,\nu)\in
O(\mu^0,W_1^{\prime},...,W_{m^{\prime}}^{\prime},\delta)\times
O(\nu^0,W_1^{\prime\prime},...W_{m^{\prime\prime}}^{\prime\prime},\delta)$
implies that
\[\min\{\frac{\mu(W^{\prime}_q)}{\mu^0(W^{\prime}_q)},
\frac{\nu(W^{\prime\prime}_s)}{\nu^0(W^{\prime\prime}_s)}
\alpha^{\prime\prime}_{qs}\}-1>-\frac{\delta}{\min\{\mu^0(W^{\prime}_q),
\nu^0(W^{\prime\prime}_s)\}}
\]and therefore
\[\lambda(W_{qs})-\lambda^0(W_{qs})>-\frac{\delta\lambda^0(W_{qs})}
{\min\{\mu^0(W^{\prime}_q),
\nu^0(W^{\prime\prime}_s)\}}>-\delta>-\varepsilon.
\]Thus, $\lambda\in O(\lambda^0,W_{11},...,W_{m^{\prime},m^{\prime\prime}},
\varepsilon)\subset O(\lambda^0,V_1,...,V_n, \varepsilon)$. The
theorem is proved.
\end{proof}

Analogically, we define the map
$\chi_{\beta}=\hat{\chi}|_{P_{\beta}(X\times Y)}$.

\begin{cor}\label{c:1tych} The map $\chi_{\beta}$ is open.
\end{cor}

The proof of the corollary is directly implied from Theorem
\ref{t:maintych} and Lemma \ref{l:2} ii).

\end{document}